\documentclass[12pt,reqno]{amsart}
\usepackage{amssymb,amsthm,amsfonts,amsmath,amscd,amssymb,epsfig,psfrag}
\usepackage{calrsfs}
\usepackage[all]{xy}

\usepackage{enumerate,array,hyperref, graphicx,color,amsthm, bbm,mathrsfs}                      
\usepackage{a4wide}
\usepackage{marginnote}           

\theoremstyle{plain}

\newtheorem{fed}{\textbf{Definition}}[section]
\newtheorem{thm}[fed]{\textbf{Theorem}}
\newtheorem{lemma}[fed]{\textbf{Lemma}}

\newtheorem{rem}[fed]{\textbf{Remark}}
\newtheorem{prop}[fed]{\textbf{Proposition}}

\newtheorem{claim}[fed]{\textbf{Claim}}
\newtheorem*{example*}{Example}

\DeclareMathAlphabet{\pazocal}{OMS}{zplm}{m}{n}

\newcommand{\N}{\mathbb{N}}

\newcommand{\Z}{\mathbb{Z}}
\newcommand{\R}{\mathbb{R}}
\newcommand{\C}{\mathbb{C}}

\newcommand{\End}{\mathrm{End}}
\newcommand{\diam}{\hspace*{\fill} $\Diamond$}

\newcommand{\vertiii}[1]{{\left\vert\kern-0.25ex\left\vert\kern-0.25ex\left\vert #1 
    \right\vert\kern-0.25ex\right\vert\kern-0.25ex\right\vert}}

\def\cA{{\pazocal A}}

\def\cF{{\pazocal F}}

\def\cP{{\pazocal P}}

\def\cV{{\pazocal V}}

\def\cX{{\pazocal X}}

\newcommand{\proofend}{\hspace*{\fill} $\Box$\\}
\def\proof{\noindent {\it Proof. \;}}

\def\1{\:\!}
\def\2{\;\!}
\def\ni{\noindent}
\def\s{\smallskip}
\def\m{\medskip}
\def\b{\bigskip}
\def\pp{\partial}
\def\id{\operatorname{id}}
\def\GL{\operatorname{GL}}
\def\Hess{\operatorname{Hess}}

\begin{document}

\title[]{A compactness result for non-local unregularized gradient flow lines}

\author{Peter Albers}
\address{Peter Albers, Mathematisches Institut, Universit\"at Heidelberg}
\email{palbers@mathi.uni-heidelberg.de}

\author{Urs Frauenfelder}
\address{Urs Frauenfelder, Mathematisches Institut, Universit\"at Augsburg}
\email{urs.frauenfelder@math.uni-augsburg.de}

\author{Felix Schlenk}  
\address{Felix Schlenk, Institut de Math\'ematiques, Universit\'e de Neuch\^atel}
\email{schlenk@unine.ch}


\date{\today}

\begin{abstract}
We prove an abstract compactness result for gradient flow lines of a 
non-local  unregularized gradient flow equation on a scale Hilbert space. 
This is the first step towards Floer theory on scale Hilbert spaces.
\end{abstract}

\maketitle
\tableofcontents

\section{Introduction}

In this article we provide the first step in the construction of non-local Floer homologies. 
Applying techniques from interpolation theory we prove compactness results 
for the space of solutions of unregularized gradient flow equations which do not need to be local, 
i.e., do not need to be solutions of a~PDE. 
Our compactness results (Theorems~\ref{th} and~\ref{strong}) are stated in the next section and proven in Sections~\ref{s:twolemmas}
and~\ref{s:proofs}.
In Section~\ref{s:floer} we show that classical Floer theory and Floer theory with delay 
on~$(\R^{2n},\omega_0)$ fit into our framework.
In the rest of this introduction we give a few motivations for why one may care 
about non-local Floer homologies. 

\subsection*{Hamiltonian Delay equations.}
A delay equation is a differential equation in which the velocity does not only depend on the present state 
but also on states in the past. 
Such equations naturally arise in population dynamics, epidemiology and economics,
but also in several problems in mechanical engineering, fluid dynamics, and visco-elasticity, 
see e.g.~\cite{Erneux09}.
In classical mechanics, delay equations arise in the (controversial) 
modified Newtonian dynamics (MOND) proposed by Milgrom~\cite{Mi83},
that serves as an alternative to the hypothesis of dark matter 
for explaining several discrepancies between observations
and theoretical computations in the dynamics of galaxies:
Under the hypothesis of Galilei invariance of MOND, 
Newton's equation becomes a delay equation, see~\cite{Mi94} and~\cite[\S 6.1.1]{FaMc12}.
%

\s
In our note~\cite{albers-frauenfelder-schlenk-delay} we addressed the question 
what a Hamiltonian delay equation may be.
The simplest Hamiltonian delay equation on~$(\R^{2n},\omega_0)$ is of the form
$$
\dot x(t) \,=\, X_H(x(t-\tau))
$$ 
where $X_H$ is the Hamiltonian vector field of 
a function $H \colon \R^{2n} \to \R$ 
and $\tau >0$ is the delay time.
On a general symplectic manifold $(M,\omega)$ with Hamiltonian function 
$H \colon M \to \R$, such an equation does not make sense, 
since $\dot x (t) \in T_{x(t)}M$ while $X_H (x (t-\tau)) \in T_{x(t-\tau)}M$. 
A general concept of a Hamiltonian delay equation on a symplectic manifold $(M,\omega)$ different from $(\R^{2n},\omega_0)$ does not exist so far.  However, by inserting into the action functional of classical mechanics a delay term, and by computing the critical point equation on the closed loop space of~$M$, we obtained many examples of delay equations that certainly deserve the predicate of Hamiltonian delay equations. 
A particular example are  delayed Lotka--Volterra equations. The search for periodic solutions of these equations with delay is an old problem in population dynamics that was already discussed in Volterra's seminal book~\cite{volterra}. 
Finding a periodic solution in this problem is somewhat relieving, 
since then there is at least one scenario in which both species survive.

Our variational approach to periodic solutions of Hamiltonian delay equations 
leads to the question whether Arnold's conjecture
on the number of periodic solutions of Hamiltonian systems 
(see for instance~\cite[Chapter 11]{mcduff-salamon17}) 
continues to hold in the delayed case. 
In particular, this would imply that the number of solutions which do not end in extinction of a species can be estimated from below by the cup-length or by the sum of the Betti numbers of the given symplectic manifold.
We in fact proved in \cite{albers-frauenfelder-schlenk-iterated-graph} that for a special 
class of delay equations, namely those which can be obtained by an iterated graph construction, the Arnold conjecture continues to hold.
This theorem only requires classical Floer theory with Lagrangian boundary conditions. 
It is therefore no question that the Arnold conjecture can be generalized to delay equations, 
but the question is how far. Indeed, one can easily cook up Hamiltonian delay equations for which 
the construction of a Floer homology is extremely doubtful. For example, if the propagation also 
depends on {\it derivatives}\/ of the solution in the past, then the Hamiltonian term in Floer's equation
is not anymore of lower order. 
But how do we actually specify what ``lower order'' means?

\subsection*{A general framework for Floer homology.}
To address the question about the meaning of ``lower order'', it is useful to briefly review 
the history of Floer homology. In their celebrated work~\cite{conley-zehnder1, conley-zehnder2} 
Conley and Zehnder proved the Arnold conjecture for the standard torus. 
With hindsight, \cite{abbondandolo-majer}, one can say that they constructed an infinite dimensional 
Morse homology on the Hilbert manifold of loops of class~$W^{1/2,2}$ on the torus. 
The result of Conley and Zehnder was later generalized by Floer~~\cite{floer} 
to a larger class of symplectic manifolds. 
Floer considered a much weaker metric for building a ``gradient vector field''
than Conley and Zehnder, namely an $L^2$-metric instead of a $W^{1/2,2}$-metric. 
In contrast to Conley and Zehnder's, Floer's ``gradient flow'' equation is not 
an ODE on a Hilbert space. Instead, by thinking of a path of loops as a cylinder,
Floer interpreted his gradient flow equation as a PDE on the cylinder in the symplectic manifold, namely a perturbed holomorphic curve equation.
Floer referred to his PDE as an ``unregularized gradient flow equation''. 

Although we nowadays have many examples of Floer homologies, see for instance~\cite{AbSch17}, 
it is still a difficult issue to say precisely what an unregularized gradient flow equation is. 
New insight into this question comes from the recent discovery 
by Hofer, Wysocki, and Zehnder~\cite{hofer-wysocki-zehnder}
of new smooth structures in infinite dimensions. 
While scales of spaces are an old topic in interpolation theory, 
see for example~\cite{triebel}, a completely unexpected result 
in~\cite{hofer-wysocki-zehnder} shows that one can define 
on a scale of spaces a new notion of smoothness, which leads to new types 
of infinite dimensional manifolds, called scale-manifolds. As with Hilbert manifolds, in the case of finite dimensions this new notion restricts to the usual notion of a manifold. However, in contrast to the setting of Hilbert manifolds, the structure of scales allows one to define a notion of ``lower order''. 
These are the $sc^+$-vector fields defined by Hofer, Wysocki, and Zehnder. 

Since Floer's equation is a PDE in the symplectic manifold, it is nowadays perceived that locality is crucial in the construction of Floer homology. 
This perception is supported by the fact that compactness
properties of the moduli space of Floer's equation are a consequence of Gromov's compactness theorem for $J$-holomorphic curves in symplectic manifolds.
In contrast, 
our guiding principle is that a Floer homology should be a Morse homology on a scale manifold. 
The generators of such a Morse complex do not need to be local solutions of an~ODE, but may contain delay or be non-local in an even stronger sense,
and the gradient flow equation is non-local a fortiori. 
A key question is if in such a setting one can still expect compactness of the solution space of the gradient flow equation. This is the issue addressed in this article, 
in which compactness properties of a certain non-local~ODE   
on a scale space are proved. 
This analysis owes much to the work of Robbin and Salamon~\cite{robbin-salamon1}
and can be thought of as a non-linear generalization of their results. 
A crucial ingredient here is the compactness of the embeddings between ascending scales. 
This compactness requirement in some way plays the role of the local compactness of 
finite-dimensional manifolds used in the construction of Morse homology 
and is a feature which is missing in the Hilbert manifold setting. 

\subsection*{Conceptual and technical advantages}
If one wishes to see the Floer homology of a symplectic manifold $(M,\omega)$ 
as the exact analogue of Morse homology on a finite-dimensional manifold, 
one should define the $L^2$-Riemannian metric on the loop space $\Lambda M$ 
in terms of an $\omega$-compatible almost complex structure that lives on $\Lambda M$, not just on~$M$, 
see the preface to the appendix of the arXiv version arXiv:1312.5201
of~\cite{FS16}.
Such an almost complex structure leads to a non-local gradient flow equation, that allows 
a much bigger space of perturbations. 
This is for example of interest if one wants to achieve transversality. 
For instance, this can be used in equivariant symplectic homology. 
The foundational and quite difficult work of Bourgeois and Oancea~\cite{bourgeois-oancea1, bourgeois-oancea2} 
shows that one can construct equivariant symplectic homology using just local methods, 
see \cite[Proposition~3.6]{bourgeois-oancea1} and 
\cite[Example~2.4]{bourgeois-oancea2}. 
On the other hand, if one allows non-local methods, the construction becomes analogous to 
the construction of equivariant Morse homology on finite-dimensional manifolds and 
can therefore easily be adapted to different equivariant Floer theoretic set-ups. 
%
%
%

Non-local gradient flow lines arise even if the critical point equation is local. 
An example is Rabinowitz--Floer homology~\cite{albers-frauenfelder, cieliebak-frauenfelder}.
Here, the critical point equation  
$$
\left\{\begin{array}{rcl}
J (v(t)) \bigl[ \pp_t v(t) - \eta \2 X_H (v(t)) \bigr] &=& 0 \\ [0.3em]
          - \int_{S^1} H \bigl (v(t) \bigr) \,dt &=& 0
\end{array}\right.
$$
for a loop $v \colon S^1 \to M$ and the Lagrangian multiplier $\eta \in \R$ is local, 
since the second seemingly non-local equation reduces to $H(v(t))=0$ for all $t\in S^1$,
but the $L^2$ gradient flow equation 
$$
\left\{\begin{array}{rcl}
\pp_s v(s,t) + J(v(s,t)) \bigl[ \pp_t v(s,t) - \eta(s) \2 X_H (v(s,t)) \bigr] &=& 0 \\ [0.3em]
\pp_s \eta (s) - \int_{S^1} H \bigl (v(s,t) \bigr) \,dt &=& 0
\end{array}\right.
$$
fails to be local:
In the second equation, $\pp_s \eta (s)$ depends on the whole loop $v(s, \cdot)$. Further, for the perturbed Rabinowitz action functional, even the critical point equation becomes non-local.

\subsubsection*{Acknowledgment}
We thank Irene Seifert very much for carefully reading the preprint. FS cordially thanks Augsburg University for its warm hospitality in the autumn of 2017.
PA is supported by DFG CRC/TRR 191, UF is supported by DFG FR/2637/2-1, and FS is supported by SNF grant 200020-144432/1.

\section{The compactness result}
In this section we formulate our compactness results.
Assume that $f \colon \mathbb{N} \to (0,\infty)$ is a monotone increasing unbounded function. 
Define the Hilbert space $\ell^2_f$ as the vector space of all real sequences 
$x=\{x_\nu\}_{\nu \in \mathbb{N}}$ satisfying
$$
\sum_{\nu=1}^\infty f(\nu) \, x_\nu^2 \,<\, \infty
$$
endowed with the inner product
$$
\langle x, y \rangle_f \,=\, \sum_{\nu=1}^\infty f(\nu) \, x_\nu y_\nu.
$$
For $k \in \mathbb{Z}$ abbreviate
$$
H_k \,=\, \ell^2_{f^k}.
$$
Note that $H_k$ is dual to $H_{-k}$ with respect to the standard inner product on $\ell^2=H_0$.
The assumption that $f$ is unbounded implies that the inclusion $H_{k+1} \to H_k$ 
is dense and compact. Let $H= \bigcup_{k \in \mathbb{Z}} H_k$. 
We choose a map $\zeta \colon \mathbb{N} \to \{ \pm 1 \}$ and define
$$
\cF \equiv \cF_\zeta \colon H \to H, \quad 
x \mapsto \Bigl\{ \zeta(\nu) \sqrt{f(\nu)} \, x_\nu \Bigr\}_{\nu \in \mathbb{N}}.
$$
For every $k \in \mathbb{Z}$ the map $\cF$ restricts to an isometry
$$
\cF \colon H_{k+1} \to H_k.
$$
We refer to $\cF$ as the \emph{fundamental operator}. The fundamental operator
can be interpreted as the Hessian of the smooth quadratic functional
$$
\mathcal{A} \colon H_1 \to \mathbb{R}, \quad 
x \mapsto \tfrac{1}{2} \sum_{\nu=1}^\infty \zeta(\nu) \sqrt{f(\nu)} \, x_\nu^2
$$
with respect to the standard inner product on $H_0=\ell^2$. 

\m
We first introduce the notion of a moving frame. 
\begin{fed}
A \emph{moving frame} is a map
$$
\Phi \colon H_1 \to \mathcal{L}(H_0,H_0)
$$
such that for every $x \in H_1$ there exists a continuous bilinear map
$$
D \Phi(x) \colon H_0 \times H_0 \to H_0
$$
such that the following properties hold.

\m
\begin{itemize}
\item[\bf $(\Phi 1)$] 
For every $x \in H_1$ the continuous linear map 
$$\Phi(x) \colon H_0 \to H_0$$
is an isomorphism.

\m
\item[\bf $(\Phi 2)$] 
The map
$$H_1 \times H_0 \to H_0, \quad (x,v) \mapsto \Phi(x) \2 v$$
is continuous.

\m
\item[\bf $(\Phi 3)$] 
For $x \in H_1$, $h \in H_1$ and $v \in H_0$, 
$$
\lim_{\|h\|_1 \to 0} \frac{1}{\|h\|_1} \cdot \| \Phi(x+h)v-\Phi(x)v-D \Phi(x)(h,v) \|_0=0.
$$

\m
\item[\bf $(\Phi 4)$] 
The map 
$$H_1 \times H_0 \times H_0 \to H_0, \quad (x,h,v) \mapsto D \Phi(x)(h,v)$$
is continuous.

\m
\item[\bf $(\Phi 5)$] 
For every $x \in H_1$ the isomorphism $\Phi(x) \colon H_0 \to H_0$ restricts to
an isomorphism $\Phi(x) \colon H_1 \to H_1$ such that the maps
$$
H_1 \times H_1 \to H_1, \quad (x,v) \mapsto \Phi(x) \2 v, \quad (x,v) \mapsto \Phi(x)^{-1} \2 v
$$
are continuous. 

\m
\item[\bf $(\Phi 6)$] 
For every $\kappa>0$ there exists a constant $c_0=c_0(\kappa)$ such that for every $x$ in the ball 
$\left\{ x \in H_1 : \|x\|_1 \leq \kappa \right\}$:
$$
\|\Phi(x)\|_{\mathcal{L}(H_0,H_0) \cap \mathcal{L}(H_1,H_1)} \leq c_0, \quad
\|\Phi(x)^{-1}\|_{\mathcal{L}(H_0,H_0) \cap \mathcal{L}(H_1,H_1)} \leq c_0
$$ 
and
$$
\|D \Phi(x)\|_{\mathcal{B}(H_0)} \leq c_0.
$$
\end{itemize}
\end{fed}

In $(\Phi 6)$ and throughout, the intersection $X \cap Y$ of two Banach spaces 
$(X, \| \: \|_X)$ and $(Y, \| \: \|_Y)$ 
is endowed with the Banach norm $\max \{ \| \: \|_X,\2 \| \: \|_Y \}$. Moreover, $\|D \Phi(x)\|_{\mathcal{B}(H_0)}$ denotes the norm as a bilinear form.

\begin{example*}
{\rm 
The trivial frame $\Phi \equiv \id_{H_0}$, $D \Phi \equiv 0$ is a moving frame.
}
\end{example*}

\begin{fed}
An \emph{unregularized vector field} $\cV \colon H_1 \to H_0$ is a map
$$
\cV \in C^0(H_1,H_0) \cap C^0(H_2,H_1)
$$
such that at every $x \in H_1$ the map $\cV \colon H_1 \to H_0$
is differentiable with differential
$$
D \cV(x) \colon H_1 \to H_0
$$
with the following properties.

\m
\begin{itemize}
\item[\bf $(\cV 1)$] 
$D \cV$ is continuous in the compact open topology, i.e., the map
$$
H_1 \times H_1 \to H_0, \quad (x,\hat x) \mapsto D \cV(x) \2 \hat x
$$
is continuous. 

\m
\item[\bf $(\cV 2)$] 
There exists a moving frame $\Phi$ such that for every $x \in H_1$ the map
$$
\Phi(x) \, D \cV(x) \, \Phi(x)^{-1} -\cF  \colon H_1 \to H_0
$$
extends to a continuous linear operator 
$$
\cP(x) \colon H_0 \to H_0
$$
with the property that the map
$$
H_1 \times H_0 \to H_0, \quad (x,\hat x) \mapsto \cP(x) \2 \hat x
$$
is continuous.

\m
\item[\bf $(\cV 3)$] 
For every $\kappa >0$ there exists a constant $c_1 = c_1(\kappa)>0$ such that 
for every $x$ in the ball 
$\left\{ x \in H_1 : \|x\|_1 \leq \kappa \right\}$:
$$
\|\cP(x)\|_{\mathcal{L}(H_0,H_0)} \leq c_1, 
\quad \mbox{and if in addition $x \in H_2$:} \;\;
\|x\|_2 \leq c_1 \bigl( \|\cV(x)\|_1 +1 \bigr) .
$$
\end{itemize}
\end{fed}
\begin{fed}
An unregularized vector field $\cV$ is called \emph{elementary} if the moving frame in assumption $(\cV 2)$ can be chosen to be the identity,
and if assumption $(\cV 3)$ can be strengthened to the assumption
\begin{itemize}
\item[\bf $(\cV 3')$] 
There exists a uniform constant $c_1' >0$ such that 
for every $x \in H_1$:
$$
\|\cP(x)\|_{\mathcal{L}(H_0,H_0)} \leq c_1'  \quad \mbox{and} \quad
  \|x\|_1 \leq c_1' \bigl( \|\cV(x)\|_0+1 \bigr),
$$
and if in addition $x \in H_2$: \;\;
$\|x\|_2 \leq c_1' \bigl( \|\cV(x)\|_1 + \|x\|_1 +1 \bigr)$.
\end{itemize}
\end{fed}

\m
In the case of the classical Floer equation on $(\R^{2n},\omega_0)$
or on the standard torus~$(T^{2n},\omega_0)$ and for
the constant almost complex structure~$i$,
the unregularized vector field is
	$\cV(x) = -i \pp_t x - \nabla H_t(x)$, 
and (roughly) 
$\cF  = -i \pp_t $ and $\cP(x) = \Hess H_t(x)$.
In our setting, the Hamiltonian term of the Floer equation, however, 
does not need to be of the classical form $-\nabla H_t(x)$,
but can be non-local, as is the case if it contains delay, 
see \S \ref{ss:delay}.
The unregularized vector field~$\cV$ above is elementary.
A moving frame $\Phi$ arises if the almost complex structure~$J$
is not the constant complex structure~$i$, 
as it happens if one writes the Floer equation on a symplectic manifold
in a symplectic chart.

\m
We fix an unregularized vector field~$\cV$ and for $T>0$ look at solutions
$w \colon I_T \to H_1$ of the equation
\begin{equation}\label{eq}
\partial_s w=\cV(w)
\end{equation}
where $I_T := (-T,T)$.

\begin{thm} \label{th}
Suppose that 
$$
w_\nu \in C^0(I_T,H_1) \cap C^1(I_T,H_0)
$$
for $\nu \in \mathbb{N}$ is a sequence of solutions of~\eqref{eq}, 
for which there exists a constant~$\kappa$ such that
\begin{equation} \label{e:ass1}
\|w_\nu\|_{C^0(I_T,H_1) \cap C^1(I_T,H_0)} \,\leq\, \kappa, \quad \forall\,\,\nu \in \mathbb{N}.
\end{equation}
Then a subsequence of $w_\nu$ converges to a solution of~\eqref{eq} in the Banach space
$C^0(I_T,H_1) \cap C^1(I_T,H_0)$.
\end{thm}

In the case that $\cV$ is elementary we get a stronger result.

\begin{thm} \label{strong}
Suppose that 
\begin{equation*} 
w_\nu \in C^0(I_T,H_1) \cap C^1(I_T,H_0)
\end{equation*}
for $\nu \in \mathbb{N}$ is a sequence of solutions of~\eqref{eq}
for an elementary unregularized vector field~$\cV$, for which there exists a 
constant~$\kappa$ such that
\begin{equation} \label{e:ass2}
\|\partial_s w_\nu\|_{L^2(I_T,H_0)} \leq \kappa, \quad \forall\,\,\nu \in \mathbb{N}.
\end{equation}
Then a subsequence of $w_\nu$ converges to a solution of~\eqref{eq} in the Banach space
$C^0(I_T,H_1) \cap C^1(I_T,H_0)$.
\end{thm}

\begin{rem} \label{re:bound}
{\rm
(i)
If in Theorem~\ref{strong} the vector field~$\cV$ is the negative $H_0$-gradient of a functional~$\cA$, 
then assumption~\eqref{e:ass2} can be guaranteed by looking only at
trajectories whose $\cA$-values (``actions'') are in a fixed compact interval. 
Indeed, if $\cV (w) = -\nabla \cA (w)$, then along a solution~$w$,
$$
\cA (w(-T)) - \cA (w(T)) \,=\,
-\int_{-T}^T \pp_s \cA (w(s)) \, ds \,=\, 
\int_{-T}^T \langle -\nabla \cA (w), \pp_s w \rangle \,ds \,=\,
\| \pp_s w \|^2_{L^2(I_T,H_0)}.
$$
An example in which this assumption is met is the Floer equation with or without delay 
on~$(\R^{2n}, \omega_0)$, 
see \S~\ref{s:floer}.

\s
(ii)
The verification of assumption~\eqref{e:ass1} is a more severe problem.
Already in classical Floer homology it is known that 
compactness of the space of gradient Floer lines cannot be
achieved by action bounds alone, due to the phenomenon of bubbling. 
On manifolds~$(M,\omega)$ with $[\omega] |_{\pi_2(M)} =0$ bubbling can be excluded, 
and it is then well-known in classical Floer homology
that on trajectories with actions in a fixed compact interval the bound~\eqref{e:ass1} holds. 
We expect that for these symplectic manifolds, the bound~\eqref{e:ass1}
can be proven also in the case of delay.
This is an interesting research project for the future. 
}
\end{rem}

\section{Two lemmas} \label{s:twolemmas}

For the proof of Theorems~\ref{th} and~\ref{strong} we need two auxiliary lemmas. 

\begin{lemma} \label{l1}
Suppose that $w \in C^0(I_T,H_1) \cap C^1(I_T,H_0)$ is a solution of~\eqref{eq}
such that there exists a constant $\kappa$ with the property that
\begin{equation}\label{est}
\|w\|_{C^0(I_T,H_1) \cap C^1(I_T,H_0)} \,\leq\, \kappa.
\end{equation}
Then  for every $0<T'<T$ it holds that $w \in \bigcap_{k=0}^2 W^{k,2}(I_{T'},H_{2-k})$
and there exists a constant $c=c(\kappa,T')$ such that
$$
\|w\|_{\bigcap_{k=0}^2 W^{k,2}(I_{T'},H_{2-k})} \,\leq\, c.
$$
\end{lemma}

\proof
Let $\Phi$ be the moving frame for $\cV$ as
in assumption~$(\cV 2)$ of the definition of an unregularized vector field~$\cV$. 
We set
\begin{equation}\label{xi}
\xi \,:=\, \Phi(w) \partial_s w = \Phi(w) \cV(w) \,\in\, C^0(I_T,H_0).
\end{equation}
From \eqref{est} and property $(\Phi 6)$ of a moving frame we observe that
\begin{equation}\label{esta}
\|\xi\|_{C^0(I_T,H_0)} \,\leq\, c_0 \kappa.
\end{equation}

\begin{claim} \label{claim}
For every $0<T'<T$ it holds that $\xi \in C^1(I_{T'},H_{-1})$ and
\begin{eqnarray} 
\label{equi}
\nonumber
\partial_s \xi 
&=& D\Phi(w) \bigl( \Phi (w)^{-1} \xi,\Phi (w)^{-1} \xi \bigr) +
\cP(w)\xi+\cF \xi .
\end{eqnarray}
\end{claim}

\proof
We fix $T' \in (0,T)$ and abbreviate $\varepsilon := \tfrac{T-T'}{2}>0$. 
Choose a smooth cutoff function $\beta \in C^\infty(I_T,[0,1])$
satisfying
$$
\beta(s) = 
\left\{\begin{array}{cl}
1, & s \in \overline{I_{T'}} \\
0, & s \in I_T \setminus I_{T-\varepsilon} .
\end{array}\right.$$
Pick a bump function $\rho \in C^\infty(\mathbb{R},[0,\infty))$ 
with the properties  
$$
\rho(\sigma)=0 \,\mbox{ for } |\sigma| \geq 1 \quad \mbox{ and } \quad 
\int_{-1}^1 \rho (\sigma) \,d \sigma =1.
$$
For $\delta>0$ we set
$\rho_\delta (s) = \tfrac{1}{\delta} \rho \bigl( \tfrac{s}{\delta} \bigr)$ and abbreviate
$$
w_\nu \,:=\,  \rho_{\varepsilon /\nu} * (\beta \2 w).
$$
Then $w_\nu$ has compact support in $I_T$, and $w_\nu \in C^\infty (I_T, H_1)$.
Now set 
$$
\xi_\nu (s) \,:=\, \Phi (w_\nu(s)) \2 \cV (w_\nu(s)) \,\in\, C^0(I_{T}, H_0).
$$
Since $\pp_s w_\nu(s) \in H_1$ for each $s\in \R$, it follows from~$(\Phi 2)$--$(\Phi 4)$ and~$(\cV 1)$
that the curve $\xi_\nu$ is differentiable in~$H_0$, and that
$$
\partial_s \xi_\nu \,=\, D\Phi(w_\nu)(\partial_s w_\nu, \cV(w_\nu)) + \Phi(w_\nu) D \cV(w_\nu) \partial_s w_\nu .
$$
This with $(\Phi 4)$ and $(\Phi 2), (\cV 1)$ imply that $\pp_s \xi_\nu \in C^0(I_{T}, H_0)$.
Using 
$$
\Phi (x) D\cV(x) \Phi(x)^{-1} \,=\, \cP(x) + \cF(x) \quad \mbox{ for $x \in H_1$}
$$ 
we can write 
\begin{eqnarray*}
\partial_s \xi_\nu &=& 
D\Phi(w_\nu)(\partial_s w_\nu, \cV(w_\nu)) + \Phi(w_\nu) D \cV(w_\nu) \partial_s w_\nu \\
&=& D\Phi(w_\nu)(\partial_s w_\nu, \cV(w_\nu)) + 
                     \cP(w_\nu) \Phi (w_\nu) \partial_s w_\nu + \cF \Phi (w_\nu) \partial_s w_\nu .
\end{eqnarray*}
We next show that in $C^0(I_{T},H_{-1})$ 
as $\nu \to \infty$ we have the convergence,
$$
\begin{array}{ccccc} 
D\Phi(w_\nu)(\partial_s w_\nu, \cV(w_\nu)) &+& \cP(w_\nu) \Phi (w_\nu) \partial_s w_\nu &+& \cF \Phi (w_\nu) \partial_s w_\nu \\
\downarrow &&\downarrow&&\downarrow \\
D\Phi(w)(\partial_s w, \cV(w)) &+& \cP(w) \Phi (w) \partial_s w &+& \cF \Phi (w) \partial_s w .
\end{array}
$$
Since on the smaller interval $I_{T'}$ we have $w_\nu \to w$ in $C^0(I_{T'}, H_1) \cap C^1 (I_{T'}, H_0)$ 
it follows that  $\cV (w_\nu) \to \cV (w)$ in $C^0(I_{T'},H_0)$
and $\Phi (w_\nu) \partial_s w_\nu \to \Phi (w) \partial_s w$ in $C^0(I_{T'},H_0)$ by $(\Phi 2)$.
Hence the left $\downarrow$ follows together with $(\Phi 4)$ in $C^0(I_{T'},H_0)$,
the middle $\downarrow$ follows together with $(\cV 2)$ in $C^0(I_{T'},H_0)$,
and the right $\downarrow$ follows in $C^0(I_{T'},H_{-1})$ since $\cF \colon H_0 \to H_{-1}$ is an isometric isomorphism.

Also note that $\xi_\nu = \Phi (w_\nu) \1 \cV (w_\nu) \to \Phi (w) \1 \cV (w) = \xi$ in 
$C^0(I_{T'}, H_0) \subset C^0(I_{T'},H_{-1})$.
We denote by 
$$
\eta \,:=\, D\Phi(w)(\partial_s w, \cV(w)) + \cP(w) \Phi (w) \partial_s w + \cF \Phi (w) \partial_s w 
$$
our candidate for $\pp_s \xi$. 
We then have
$\xi_\nu \to \xi$ and $\pp_s \xi_\nu \to \eta$ in $C^0(I_{T'},H_{-1})$.
Therefore, with $\langle \:,\2 \rangle$ the $H_0$ inner product, 
\begin{eqnarray*}
\int_{I_{T'}} \langle \eta, \varphi \rangle \,ds \,=\,
\lim_{\nu \to \infty} \int_{I_{T'}} \langle \pp_s \xi_\nu, \varphi \rangle \,ds \,=\,
-\lim_{\nu \to \infty} \int_{I_{T'}} \langle \xi_\nu, \pp_s \varphi \rangle \,ds \,=\,
- \int_{I_{T'}} \langle \xi, \pp_s \varphi \rangle \,ds 
\end{eqnarray*}
for every compactly supported test function $\varphi \in C^\infty (I_{T'}, H_1)$ to the dual space $H_1$
of~$H_{-1}$.
Hence $\eta \in C^0(I_{T'},H_{-1})$ is a weak derivative of $\xi$ and so, being continuous, 
is its derivative.
\proofend

\begin{claim}\label{claim1.5}
It holds that $\xi \in L^2(I_T,H_1) \cap W^{1,2}(I_T,H_0)$
and we have the estimates
\begin{eqnarray} \label{est9}
\|\partial_s \xi\|_{L^2(I_{T'},H_0)} &\leq&  
                 c_0 \kappa \Bigl( \|\partial_s \beta\|_2 + \sqrt{2T} (\kappa + c_1) \Bigr) , 
\\
\label{est10}
\|\xi\|_{L^2(I_{T'},H_1)} &\leq& 
  c_0 \kappa \Bigl(\|\partial_s \beta\|_2 + \sqrt{2T} (\kappa + c_1) \Bigr). 
\end{eqnarray}
\end{claim}

\proof
We define $\beta$, $\rho$, $\rho_\delta$ as before, and for $0 < \delta < \tfrac{\varepsilon}{2}$ set
$\xi^\beta = \beta \2 \xi$ and
$$
\xi^\beta_\delta:=\rho_\delta * \xi^\beta \,\in\, C^\infty (I_T, H_0).
$$
Then $\xi^\beta_\delta$ has compact support in $I_T$.
From Claim~\ref{claim} we obtain
$$
\xi \,=\, \cF^{-1} \partial_s \xi-\cF^{-1}D \Phi(w)
\bigl( \Phi(w)^{-1}\xi,\Phi(w)^{-1}\xi \bigr) -\cF^{-1} \cP(w) \xi .
$$
Therefore
\begin{eqnarray*}
\xi_\delta^\beta &=& 
\rho_\delta *\Bigl(\cF^{-1} \beta \partial_s \xi-\cF^{-1}D \Phi(w)
\bigl(\Phi(w)^{-1}\xi,\Phi(w)^{-1}\xi^\beta\bigr)
-\cF^{-1} \cP(w) \xi^\beta \Bigr) \\
&=&(\partial_s \rho_\delta) *(\cF^{-1} \xi^\beta)
                                -\rho_\delta * \cF^{-1} \Bigl( (\partial_s \beta)\xi+R \Bigr) 
\end{eqnarray*}
where we abbreviate
$$
R \,:=\, D \Phi(w) \bigl( \Phi(w)^{-1}\xi,\Phi(w)^{-1}\xi^\beta \bigr)
 + \cP(w)\xi^\beta.
$$
Since $\cF \colon H_1 \to H_0$ and $\cF \colon H_0 \to H_{-1}$ are isometric isomorphisms, 
this formula and the properties of~$\Phi$ and~$\cP$ imply that
$$
\xi_\delta ^\beta \,\in\, C^\infty (I_T, H_1).
$$
Moreover, we compute
\begin{eqnarray} \label{est1}
\nonumber
\partial_s \xi_\delta^\beta - \cF\xi_\delta^\beta
&=&
\partial_s \xi^\beta_\delta-\cF \Bigl( (\partial_s \rho_\delta)*(\cF^{-1} \xi^\beta) \Bigr) + 
         \cF \Bigl( \rho_\delta * \cF^{-1}\Bigl( (\partial_s \beta) \xi + R \Bigr) \Bigr) \\ \nonumber
&=&
\partial_s \xi_\delta^\beta-\partial_s \xi_\delta^\beta + 
                                \rho_\delta * \Bigl( (\partial_s \beta) \xi+R \Bigr)\\ 
&=&
\rho_\delta * \Bigl((\partial_s \beta)\xi+R\Bigr). 
\end{eqnarray}
By Young's inequality,
$$
\bigl\| \rho_\delta * \bigl( (\partial_s \beta) \xi+R \bigr) \bigr\|_{L^2(I_T,H_0)}
\,\leq\,
\bigl\| \rho_\delta \bigr\|_{L^1(I_T,\mathbb{R})}
\cdot \bigl\| (\partial_s \beta) \xi+R \bigr\|_{L^2(I_T,H_0)}.
$$
Since
$$
\|\rho_\delta\|_{L^1} \,=\, \int_{-T}^T \rho_\delta(s) \, ds \,=\, 
\int_{-T}^T \tfrac{1}{\delta}\rho \bigl( \tfrac{s}{\delta} \bigr) \, ds 
\,=\, \int_{-T/\delta}^{T/\delta} \rho(\sigma) \, d \sigma \,=\, 1
$$
we can estimate using \eqref{esta}, the bound on $D\Phi (w(s))$ in~$(\Phi 6)$,
the definition $\Phi (w)^{-1}\xi = \partial_s w$, assumption~\eqref{est},
and property~$(\cV 3)$:
\begin{eqnarray}\label{est2}\nonumber
\Bigl\| \rho_\delta * \Bigl( (\partial_s \beta)\xi + R \Bigr) \Bigr\|_{L^2(I_T,H_0)}&\leq&  \Bigl\| (\partial_s \beta) \xi+R \Bigr\|_{L^2(I_T,H_0)}\\ \nonumber
&\leq& c_0 \kappa \|\partial_s \beta\|_2 + \sqrt{2T}(c_0\kappa^2+c_0 c_1 \kappa) \\ 
&=&c_0 \kappa \Bigl( \|\partial_s \beta\|_2 + \sqrt{2T} (\kappa+c_1) \Bigr).
\end{eqnarray}
From \eqref{est1} and \eqref{est2} we infer
\begin{equation} \label{est4} 
\|\partial_s \xi_\delta^\beta - \cF \xi_\delta^\beta\|_{L^2(I_T,H_0)}
\,\leq\, c_0 \kappa \Bigl( \|\partial_s \beta\|_2 + \sqrt{2T}(\kappa+c_1) \Bigr).
\end{equation}
Note that $\cF$ is selfadjoint with respect to~$\langle \cdot, \cdot \rangle$, the inner product on~$H_0$.
Since $\xi^\beta_\delta$ has compact support we obtain using integration by parts
\begin{eqnarray*}
\int_{-T}^T \langle \partial_s \xi^\beta_\delta, \cF \xi^\beta_\delta \rangle \, ds
&=&-\int_{-T}^T \langle \xi^\beta_\delta, \cF \partial_s \xi^\beta_\delta \rangle \, ds\\
&=&-\int_{-T}^T \langle \cF \xi^\beta_\delta, \partial_s \xi^\beta_\delta \rangle \, ds\\
&=&-\int_{-T}^T \langle \partial_s \xi^\beta_\delta, \cF \xi^\beta_\delta \rangle \, ds,
\end{eqnarray*}
from which we deduce
$$\int_{-T}^T \langle \partial_s \xi^\beta_\delta, \cF \xi^\beta_\delta \rangle \, ds \,=\, 0.$$
This implies that
\begin{eqnarray*}
\| \partial_s \xi^\beta_\delta - \cF \xi^\beta_\delta\|^2_{L^2(I_T,H_0)}
&=&
\int_{-T}^T \big\langle \partial_s \xi^\beta_\delta-\cF \xi^\beta_\delta,
\partial_s \xi^\beta_\delta-\cF\xi^\beta_\delta \big\rangle \, ds\\
&=&\int_{-T}^T \big \langle \partial_s \xi^\beta_\delta,\partial_s \xi^\beta_\delta \big \rangle \, ds 
- 2 \int_{-T}^T \big \langle \partial_s \xi^\beta_\delta, \cF\xi^\beta_\delta \big \rangle \, ds \\
& &+\int_{-T}^T \big \langle \cF \xi^\beta_\delta, \cF\xi^\beta_\delta \big \rangle \, ds \\
&=&\|\partial_s \xi^\beta_\delta\|_{L^2(I_T,H_0)}^2+\|\cF\xi^\beta_\delta\|^2_{L^2(I_T,H_0)}\\
&=&\|\partial_s \xi^\beta_\delta\|_{L^2(I_T,H_0)}^2+\|\xi^\beta_\delta\|^2_{L^2(I_T,H_1)}. 
\end{eqnarray*}
Here we have used for the last equality that $\cF \colon H_1 \to H_0$ is an isometric isomorphism.
Combining this with~\eqref{est4} we obtain the two estimates
\begin{eqnarray} 
\label{est5}
\|\partial_s \xi_\delta^\beta\|_{L^2(I_T,H_0)} &\leq& 
              c_0 \kappa \Bigl( \|\partial_s \beta\|_2 + \sqrt{2T}(\kappa+c_1) \Bigr) , \\
\label{est6}
\|\xi_\delta^\beta\|_{L^2(I_T,H_1)} &\leq&
                   c_0 \kappa \Bigl( \|\partial_s \beta\|_2 + \sqrt{2T} (\kappa+c_1) \Bigr).
\end{eqnarray}
Note that the bounds in \eqref{est5} and \eqref{est6} are independent of~$\delta$. 
Hence there exists a sequence $\delta_\nu \to 0$ as $\nu \to \infty$ such that  
$\xi_{\delta_\nu}^\beta$ converges weakly in $L^2(I_T,H_1) \cap W^{1,2}(I_T,H_0)$ to some
$$\xi^\beta_0 \in L^2(I_T,H_1) \cap W^{1,2}(I_T,H_0).$$
Because $\xi^\beta_{\delta_\nu}$ converges strongly in $L^2(I_T,H_0)$ to $\xi^\beta$ we conclude
that
$$\xi^\beta=\xi^\beta_0 \in L^2(I_T,H_1) \cap W^{1,2}(I_T,H_0).$$
We write $\pp_s \xi^\beta\in L^2(I_T,H_0)$ for the weak derivative of~$\xi^\beta$.
Moreover, from \eqref{est5} and~\eqref{est6} we get the estimates
\begin{eqnarray*} \label{est7}
\|\partial_s \xi^\beta\|_{L^2(I_T,H_0)} &\leq&
c_0 \kappa \Bigl( \|\partial_s \beta\|_2 + \sqrt{2T} (\kappa + c_1) \Bigr) , \\
\label{est8}
\; \|\xi^\beta\|_{L^2(I_T,H_1)} &\leq& c_0 \kappa \Bigl( \|\partial_s \beta\|_2 + \sqrt{2T} (\kappa + c_1) \Bigr).
\end{eqnarray*}
Since $\xi^\beta|_{[-T',T']}=\xi|_{[-T',T']}$, these inequalities imply Claim \ref{claim1.5}. \proofend
%

\begin{claim} \label{claim2}
For every $T'<T$ the function $\pp_s w = \Phi (w)^{-1} \xi$ has a weak derivative in $L^2(I_{T'},H_0)$,
namely
$$
\partial_s^2w\,=\,\Phi (w)^{-1} \partial_s \xi - \Phi (w)^{-1} D\Phi (w) (\pp_sw, \pp_sw) .
$$
\end{claim}
  
\proof
Abbreviate $A(s) = \Phi (w(s))$. We need to show that
$$
A(s)^{-1} \pp_s \xi - A(s)^{-1} \circ DA(s) \circ A(s)^{-1} \, \xi (s)
$$
is a weak derivative of $\pp_s w = A(s)^{-1} \xi$.
As in the proof of Claim~\ref{claim} we take $w_\nu := \rho_{\varepsilon/\nu} \ast (\beta w) \in C^\infty(I_T,H_1)$,
and we now set $\xi_\nu := \Phi (w_\nu) \pp_s w_\nu$.
With $A_\nu (s) := \Phi (w_\nu (s))$ we then have $\pp_s w_\nu = A_\nu^{-1} \xi_\nu$.
By $(\Phi 3)$, $(\Phi 4)$, the map $s \mapsto A_\nu (s) = \Phi (w_\nu (s)) \in \mathcal{L} (H_0,H_0)$
is differentiable, and so 
\begin{equation} \label{e:swA}
\pp_s^2 w_\nu \,=\, \pp_s (A_\nu^{-1} \xi_\nu) \,=\, A_\nu (s)^{-1} \pp_s \xi_\nu - A_\nu (s)^{-1} \circ D A_\nu(s) \circ A_\nu(s)^{-1} \,\xi_\nu (s) \,\in\, H_0 .
\end{equation}
In order to prove Claim \ref{claim2} we will argue below that for every compactly supported 
$\varphi \in C^\infty (I_{T'},H_0)$ it holds that
\begin{eqnarray}
\int \left\langle A(s)^{-1} \xi, \pp_s \varphi \right\rangle  &=& 
\phantom{-} \lim_{\nu \to \infty} \int \left\langle A_\nu(s)^{-1} \xi_\nu, \pp_s \varphi \right\rangle 
\label{E1}\\
&=&  
- \lim_{\nu \to \infty} \int \left\langle A_\nu(s)^{-1} \pp_s \xi_\nu - A_\nu(s)^{-1} DA_\nu(s) A_\nu(s)^{-1} \xi_\nu(s), \varphi \right\rangle \label{E2} \\
&=& \phantom{\lim_{\nu \to \infty} \!} - \int \left\langle A(s)^{-1} \pp_s \xi - A(s)^{-1} \, DA(s) \, A(s)^{-1} \xi(s), \varphi \right\rangle \label{E3}
\end{eqnarray}
where all integrals are over $I_{T'}$ and $\langle \; , \, \rangle = \langle \; , \, \rangle_{H_0}$.

\m \ni
\eqref{E1}: Recall from the proof of Claim \ref{claim} that $\Phi (w_\nu) \pp_s w_\nu \to \Phi (w) \pp_s w$
in $C^0(I_{T'},H_0)$, that is, 
$$
\xi_\nu \to \xi \; \mbox{ in } C^0(I_{T'}, H_0).
$$
Further, by $(\Phi 6)$, the functions $\left| \langle A_\nu(s)^{-1} \xi_\nu, \pp_s \varphi \rangle \right|$
are uniformly bounded by the constant $c_0 \left( \|\xi\|_{C_0(I_{T'},H_0)} +1 \right) \|\pp_s \varphi\|_{C_0(I_{T'},H_0)}$
for $\nu$ large enough. Hence \eqref{E1} follows from the dominated convergence theorem.

\m \ni
\eqref{E2}: By $(\Phi 2)$--$(\Phi 4)$,  
$$\pp_s \xi_\nu \,=\, \pp_s \bigl( \Phi (w_\nu) \pp_s w_\nu \bigr) \,=\,
D\Phi (w_\nu) (\pp_sw_\nu, \pp_s w_\nu) + \Phi (w_\nu) \pp_s^2 w_\nu \;\in\; C^0(I_{T'},H_0) .
$$
Hence \eqref{E2} follows from \eqref{e:swA}.

\m \ni
\eqref{E3}: The identity
$$
\lim_{\nu \to \infty} \int \left\langle A_\nu(s)^{-1} DA_\nu(s) A_\nu(s)^{-1} \xi_\nu(s), \varphi \right\rangle 
\,=\,
\int \left\langle A(s)^{-1} \, DA(s) \, A(s)^{-1} \xi(s), \varphi \right\rangle
$$
follows from the continuity of $\Phi^{-1}$ and $D \Phi$, from $\xi_\nu \to \xi$ in $C^0(I_{T'}, H_0)$, 
and from dominated convergence.
Further, 
since $\pp_s \xi$ is the weak $L^2(I_{T'}, H_0)$ derivative of $\xi$, the identity
$$
\lim_{\nu \to \infty} \int \left\langle A_\nu(s)^{-1} \pp_s \xi_\nu, \varphi \right\rangle 
\,=\,
\int \left\langle A(s)^{-1} \pp_s \xi, \varphi \right\rangle
$$
is equivalent to 
\begin{equation} \label{e:last}
\lim_{\nu \to \infty} \int \left\langle \xi_\nu, \pp_s \left([A_\nu(s)^{-1}]^T \varphi \right) \right\rangle 
\,=\,
\int \left\langle \xi, \pp_s \left( [A(s)^{-1}]^T \varphi \right) \right\rangle .
\end{equation}
Note that transposition $T$ commutes with differentiation and does not change the operator norm.
The bounds in~$(\Phi 6)$ thus imply uniform bounds on the operator norms of $\pp_s [A_\nu(s)^{-1}]^T$
and $[A_\nu(s)^{-1}]^T$.
Therefore, 
\eqref{e:last} again follows from the continuity of $\Phi^{-1}$ and $D \Phi$, 
from $\xi_\nu \to \xi$ in $C^0(I_{T'}, H_0)$, 
and from dominated convergence.
\proofend

Recall that we denote by $\pp_s^2 w$ the weak derivate of $\pp_s w$.
Combining Claim~\ref{claim2}
with \eqref{est}, \eqref{est9}, and with property~$(\Phi 6)$ of a moving frame,
we estimate
\begin{eqnarray*}
\|\partial_s^2 w\|_{L^2(I_{T'},H_0)}& \leq&\|\Phi(w)^{-1} \partial_s \xi\|_{L^2(I_{T'},H_0)}\\
& &+\|\Phi(w)^{-1} D \Phi(w)(\partial_s w,\partial_s w)\|_{L^2(I_{T'},H_0)}\\
&\leq&c_0\|\partial_s \xi\|_{L^2(I_{T'},H_0)}+\sqrt{2T}c_0^2 \kappa^2\\
&\leq&c_0^2 \kappa \Bigl(\|\partial_s \beta\|_2+\sqrt{2T} (c_1 + 2\kappa) \Bigr).
\end{eqnarray*} 
Combining this estimate once more with \eqref{est} we get
\begin{eqnarray} \label{main1}
\nonumber
\|w\|^2_{W^{2,2}(I_{T'},H_0)}&=&\|\partial_s^2 w\|^2_{L^2(I_{T'},H_0)}+
\|w\|^2_{W^{1,2}(I_{T'},H_0)}\\ \nonumber
&\leq&\|\partial_s^2 w\|^2_{L^2(I_{T'},H_0)}+
2T\|w\|^2_{C^1(I_T,H_0)} \\
&\leq& c_0^4 \kappa^2 \Bigl( \|\partial_s \beta\|_2+\sqrt{2T}(c_1+2\kappa) \Bigr)^2 + 2T \kappa^2 =: \kappa_0^2 .
\end{eqnarray} 
From \eqref{xi}, \eqref{est10}, and property~$(\Phi 6)$ of a moving frame we obtain
\begin{eqnarray*} 
\|\partial_s w\|_{L^2(I_{T'},H_1)}&=&\|\Phi(w)^{-1} \xi\|_{L^2(I_{T'},H_1)}\\
&\leq&c_0^2 \kappa \Bigl( \|\partial_s \beta\|_2 + \sqrt{2T} (\kappa + c_1) \Bigr).
\end{eqnarray*} 
Combining this estimate with \eqref{est} we infer 
\begin{eqnarray} \label{main2}
\nonumber
\|w\|^2_{W^{1,2}(I_{T'},H_1)} &=& \|\partial_s w\|^2_{L^2(I_{T'},H_1)}+\|w\|^2_{L^2(I_{T'},H_1)} \nonumber \\ 
&\leq& \|\partial_s w\|^2_{L^2(I_{T'},H_1)}+2T\|w\|^2_{C^0(I_T,H_1)} \nonumber\\
&\leq& c_0^4 \kappa^2 \Bigl( \|\partial_s \beta\|_2+2\sqrt{2T} (\kappa + c_1) \Bigr)^2 + 2T \kappa^2 =: \kappa_1^2 .
\end{eqnarray} 
Recalling that $\partial_s w = \cV(w)$ we estimate 
by taking advantage of assumption~$(\cV 3)$ on an unregularized vector field,
and with $\mu := \max \{ \sqrt{2T'}, 1\}$,
\begin{eqnarray} \label{main3} \nonumber
\|w\|_{L^2(I_{T'},H_2)}&\leq& c_1 \2 \mu \bigl( \|\cV(w)\|_{L^2(I_{T'},H_1)} + 1 \bigr) \\ \nonumber
&=&    c_1 \2 \mu \bigl( \|\partial_s w\|_{L^2(I_{T'},H_1)} +1 \bigr) \\ \nonumber
&\leq& c_1 \2 \mu \bigl( \|w\|_{W^{1,2}(I_{T'},H_1)} + 1 \bigr) \\ 
&\leq& c_1 \2 \mu (\kappa_1+1)=:\kappa_2.
\end{eqnarray}
Lemma \ref{l1} follows from \eqref{main1}, \eqref{main2}, and \eqref{main3} when setting $c = \max \{\kappa_0,\kappa_1,\kappa_2\}$.
\proofend

\begin{lemma} \label{l4}
For $T>0$, $p>1$, and $\ell \in \mathbb{N}$ the inclusion 
$$
\iota \colon \bigcap_{k=0}^\ell W^{k,p}(I_T,H_{\ell-k}) \,\to\, \bigcap_{k=0}^{\ell-1}C^k(I_T,H_{\ell-1-k})
$$ 
is a compact operator. 
\end{lemma}

\proof
For $N \in \mathbb{N}$ let $V_N \subset H_0=\ell^2$ be the 
$N$-dimensional subspace of sequences $x=\{x_\nu\}_{\nu \in \mathbb{N}}$ with $x_\nu=0$ for $\nu>N$. 
Observe that $V_N \subset H_k$ for every $k \in \mathbb{Z}$. 
Let 
$$\pi_N \colon H_0 \to V_N$$
be the orthogonal projection. The standard basis of $H_0=\ell^2$ is a common orthogonal basis
of~$H_k$ for every $k \in \mathbb{Z}$. In particular, the restriction
$$
\pi_N|_{H_k} \colon H_k \to V_N
$$
is the orthogonal projection of $H_k$ to $V_N$ for every $k \in \mathbb{N}$. 
Let 
$$
\Pi_N \colon \bigcap_{k=0}^\ell W^{k,p}(I_T,H_{\ell-k}) \to W^{\ell,p}(I_T,V_N), 
\quad w \mapsto \pi_N \circ w.
$$
Since $V_N$ is finite-dimensional and $p>1$, the inclusion
$$
I_N \colon W^{\ell,p}(I_T,V_N) \,\to\, C^{\ell-1}(I_T,V_N)
$$
is a compact operator. We abbreviate by
$$
J_N \colon C^{\ell-1}(I_T,V_N) \to \bigcap_{k=0}^{\ell-1} C^k(I_T,H_{\ell-1-k})
$$
the inclusion
and by 
$$
\iota_N \colon \bigcap_{k=0}^\ell W^{k,p}(I_T,H_{\ell-k}) \to \bigcap_{k=0}^{\ell-1}C^k(I_T,H_{\ell-1-k})
$$
the composition of these three maps,
$$
\iota_N \,:=\, J_N \circ I_N \circ \Pi_N.
$$
Since $I_N$ is compact and the other two maps are continuous, $\iota_N$ is a compact operator.
We are thus left with showing that $\iota_N$ converges to~$\iota$ in the norm topology as $N \to \infty$. 
Arguing by contradiction we assume that there exists a constant~$c>0$ such that 
for infinitely many $N \in \mathbb{N}$
there exists $w_N \in \bigcap_{k=0}^\ell W^{k,p}(I_T,H_{\ell-k})$ with the property that
\begin{equation} \label{wn1}
\max \Big\{ \left\| (\iota-\iota_N) w_N \right\|_{C^k(I_T,H_{\ell-1-k})} : 0 \leq k \leq \ell-1 \Big\} \,=\, 1
\end{equation}
but
\begin{equation} \label{wn2}
\max \Big\{ \| w_N\|_{W^{k,p}(I_T,H_{\ell-k})} : 0 \leq k \leq \ell \Big\} \,\leq\, c.
\end{equation}
From \eqref{wn1} we deduce that there exists $0 \leq j \leq \ell-1$
such that
$$
\left\| (\iota-\iota_N) w_N \right\|_{C^j(I_T,H_{\ell-1-j})} \,=\, 1
$$
for an infinite number of $N \in \N$.
Fix such an~$N$. Then there exists $s = s_N \in I_T$ with the property that
\begin{equation} \label{e:=1}
\left\| (\mathrm{id}-\pi_N)\partial^j_s w_N(s) \right\|_{H_{\ell-1-j}} \,=\, 1.
\end{equation}
Let $q \in (1,\infty)$ be the number dual to $p$ in the sense that 
$$\frac{1}{p}+\frac{1}{q}\,=\,1.$$
Suppose that $s' \in I_T$ satisfies $|s'-s| \leq \left( \frac{1}{2c} \right)^q$. Then with~\eqref{wn2},
\begin{eqnarray*}
|s'-s| &\leq& \bigg(\frac{1}{2c}\bigg)^q\\
&\leq& \Bigg(\frac{1}{2\|w_N\|_{W^{1+j,p}(I_T,H_{\ell-1-j})}}\Bigg)^q\\
&\leq& \Bigg(\frac{1}{2\|(\mathrm{id}-\pi_N)w_N\|_{W^{1+j,p}(I_T,H_{\ell-1-j})}}\Bigg)^q.
\end{eqnarray*}
Using also \eqref{e:=1} and H\"older's inequality, we estimate
\begin{eqnarray*}
\left\| (\mathrm{id}-\pi_N) \partial^j_s w_N(s') \right\|_{H_{\ell-1-j}}
&\geq& \left\| (\mathrm{id}-\pi_N) \partial^j_s w_N(s) \right\|_{H_{\ell-1-j}}\\
& &- \biggl| \int_s^{s'} \left\| (\mathrm{id}-\pi_N) \partial^{j+1}_s w_N(\sigma) \right\|_{H_{\ell-1-j}} \, d \sigma \biggr| \\
&\geq&1- \left\|(\mathrm{id}-\pi_N) \partial^{j+1}_s w_N \right\|_{L^p(I_T,H_{\ell-1-j})} \,
                                                                       |s'-s|^{\frac{1}{q}}\\
&\geq&1- \left\| (\mathrm{id}-\pi_N) w_N \right\|_{W^{1+j,p}(I_T,H_{\ell-1-j})} \,|s'-s|^{\frac{1}{q}}\\
&\geq&1-\tfrac{1}{2}\\
&=&\tfrac{1}{2}.
\end{eqnarray*}
Since the function $f$ is monotone increasing by assumption, we thus obtain
$$
\left\| (\mathrm{id}-\pi_N) \partial^j_s w_N(s') \right\|_{H_{\ell-j}} \,\geq\, \tfrac{1}{2}\sqrt{f(N+1)}.
$$
Using this we can estimate
\begin{eqnarray*}
\left\| w_N \right\|_{W^{j,p}(I_T,H_{\ell-j})} &\geq& 
\left\| \partial^j_s w_N \right\|_{L^p(I_T,H_{\ell-j})}  \\
&\geq& \left\| (\mathrm{id}-\pi_N) \partial^j_s w_N \right\|_{L^p(I_T,H_{\ell-j})}\\
&\geq& \tfrac 12 \sqrt{f(N+1)} \cdot \min \Big\{\bigl( \tfrac{1}{2c} \bigr)^{\frac{q}{p}},T^{\frac{1}{p}}\Big\}\\
&=&\tfrac 12 \sqrt{f(N+1)} \cdot \min \Big\{ \bigl( \tfrac{1}{2c} \bigr)^{\frac{1}{p-1}}, T^{\frac{1}{p}}\Big\}.
\end{eqnarray*}
Since $f$ is unbounded, this violates~\eqref{wn2} for $N$ large enough. 
This contradiction proves the lemma. 
\proofend

\section{Proof of Theorem~\ref{th} and \ref{strong}} \label{s:proofs}

\ni
{\it Proof of Theorem~\ref{th}.}
According to the assumptions of Theorem~\ref{th} and by Lemma~\ref{l1},
for every $0<T'<T$ the sequence $w_\nu|_{[-T',T']}$ is uniformly bounded in 
$\bigcap_{k=0}^2 W^{k,2}(I_{T'},H_{2-k})$. 
Lemma~\ref{l4} tells us that the inclusion 
$$
\bigcap_{k=0}^2 W^{k,2}(I_{T'},H_{2-k}) \,\to\, 
\bigcap_{k=0}^1 C^k(I_{T'},H_{1-k})
$$ 
is compact. Therefore, $w_\nu|_{I_{T'}}$
has a convergent subsequence in $C^0(I_{T'},H_1) \cap C^1(I_{T'}, H_0)$. 
The theorem now follows by a diagonal argument. 
\proofend

To see how Theorem~\ref{th} implies Theorem~\ref{strong} we need the following proposition.

\begin{prop}\label{pr}
Suppose that $w \in C^0(I_T,H_1) \cap C^1(I_T,H_0)$ is a solution of~\eqref{eq} for an elementary unregularized vector field $\cV$, 
such that there exists a constant $\kappa$ with the property that
$$
\|\partial_s w\|_{L^2(I_T,H_0)} \,\leq\, \kappa .
$$
Then for every $0<T'<T$ there exists a constant $c=c(\kappa,T')$ such that 
$$
\|w\|_{C^0(I_{T'},H_1) \cap C^1(I_{T'},H_0)} \,\leq\, c.
$$
\end{prop}
The proof of Proposition~\ref{pr} is based on the next two lemmas.

\begin{lemma} \label{l2}
Under the assumptions of Proposition~\ref{pr},
there exists a constant $c=c(\kappa)$ such that
$$
\|w\|_{L^2(I_T,H_1) \cap W^{1,2}(I_T,H_0)} \,\leq\, c.
$$
\end{lemma}

\proof 
By assumption, $\partial_s w = \cV(w)$, and by $(\cV 3')$, 
$$
\|w(s)\|_{H_1} \,\leq\, c_1' \left( \| \cV(w(s)) \|_{H_0} +1 \right) \quad \forall \, s \in I_T .
$$
Hence, with Cauchy--Schwarz, 
\begin{eqnarray*}
\|w\|_{L^2(I_T,H_1)}&\leq&c_1' \left(\|\cV(w)\|_{L^2(I_T,H_0)}+\sqrt{2T}\right)\\
&=&c_1' \left(\|\partial_s w\|_{L^2(I_T,H_0)}+\sqrt{2T}\right)\\
&\leq& c_1' (\kappa+\sqrt{2T}).
\end{eqnarray*}
With this and the assumption, 
\begin{eqnarray*}
\|w\|_{W^{1,2}(I_T,H_0)}&=&\sqrt{\|w\|_{L^2(I_T,H_0)}^2+\|\partial_s w\|_{L^2(I_T,H_0)}^2}\\
&\leq& \sqrt{\|w\|_{L^2(I_T,H_1)}^2+\|\partial_s w\|_{L^2(I_T,H_0)}^2}\\
&\leq&\sqrt{(c_1')^2(\kappa+\sqrt{2T})^2+\kappa^2}.
\end{eqnarray*}
Combining these two inequalities we obtain
$$
\|w\|_{L^2(I_T,H_1) \cap W^{1,2}(I_T,H_0)} \,\leq\, \sqrt{(c_1')^2(\kappa+\sqrt{2T})^2+\kappa^2}
\,=:\, c .
$$
The lemma follows.
\proofend

In the proof of Lemma~\ref{l2} we did not use the fact that the 
moving frame for $\cV$ is trivial.
This becomes crucial in the following lemma, however.

\begin{lemma}\label{l3}
Suppose that $w \in C^0(I_T,H_1) \cap C^1(I_T,H_0)$ is a solution of \eqref{eq} 
for an elementary unregularized vector field~$\cV$,
such that there exists a constant $\kappa$ with the property that
\begin{equation} \label{esti}
\|w\|_{L^2(I_T,H_1) \cap W^{1,2}(I_T,H_0)} \leq \kappa.
\end{equation}
Then for every $0<T'<T$ there exists a constant $c=c(\kappa,T')$ such that
$$
\|w\|_{\bigcap_{k=0}^2 W^{k,2}(I_{T'},H_{2-k})} \,\leq\, c.
$$
\end{lemma}

\proof
Abbreviate 
$\xi := \partial_s w = \cV(w) \in L^2(I_T,H_0)$.
By \eqref{esti},
\begin{equation}\label{esti1}
\|\xi\|_{L^2(I_T,H_0)} \leq \kappa.
\end{equation}
Since $\cV$ is elementary, its moving frame can be chosen trivial and so assumption~$(\cV 2)$ in the definition of an unregularized vector field has the simple form
$$
D\cV \,=\, \cF + \cP .
$$
Differentiating $\partial_s w = \cV (w)$ we thus find, as in Claim~\ref{claim},
\begin{equation} \label{simd}
\partial_s \xi \,=\, \cF (w) \xi + \cP (w) \xi \,\in\, L^2(I_T,H_{-1})
\end{equation}
so that 
$$
\xi \in L^2(I_T,H_0) \cap W^{1,2}(I_T,H_{-1}) .
$$
Choose $\beta$ and $\rho_\delta$ as in the proof of Lemma~\ref{l1} and introduce the
compactly supported functions
$$
\xi^\beta := \beta \2 \xi, \quad \xi_\delta^\beta := \rho_\delta * \xi^\beta
\,\in\, C^\infty (I_T, H_0) .
$$
From \eqref{simd} we see that
$$
\xi \,=\, \cF^{-1} \partial_s \xi - \cF^{-1} \cP (w) \xi
$$
so that
\begin{eqnarray*}
\xi_\delta^\beta&=&\rho_\delta * \Bigl(\cF^{-1} \beta \partial_s \xi-
\cF^{-1} \cP (w) \xi^\beta \Bigr)\\
&=&(\partial_s \rho_\delta) *(\cF^{-1} \xi^\beta)-
\rho_\delta * \cF^{-1} \Bigl( (\partial_s \beta) \xi + \cP (w) \xi^\beta \Bigr).
\end{eqnarray*}
In particular,
$$
\xi_\delta^\beta \,\in\,  C^\infty (I_T, H_1) .
$$
Identity \eqref{est1} now becomes
\begin{eqnarray} \label{sim1}
\partial_s \xi_\delta^\beta-\cF \xi_\delta^\beta
&=&\rho_\delta * \Bigl( (\partial_s \beta) \xi+ \cP (w) \xi^\beta \Bigr).
\end{eqnarray}
By Young's inequality, \eqref{esti} and property $(\cV 3')$ we can estimate
\begin{eqnarray}\label{sim2}
& &\Bigl\| \rho_\delta * \Bigl( (\partial_s \beta) \xi+ \cP (w) \xi^\beta \Bigr) \Bigr\|_{L^2(I_T,H_0)} \\ \nonumber
&=&\Bigl\|(\partial_s \beta)\xi + \cP (w)\xi^\beta\Bigr\|_{L^2(I_T,H_0)}\\ \nonumber
&\leq&\kappa\bigl( \|\partial_s \beta\|_\infty + c_1' \bigr).
\end{eqnarray}
From \eqref{sim1} and \eqref{sim2} we infer that 
$$
\| \partial_s \xi_\delta^\beta - \cF \xi_\delta^\beta\|_{L^2(I_T,H_0)}
\,\leq\, \kappa \bigl( \|\partial_s \beta\|_\infty + c_1' \bigr).
$$
In the same way as we deduced \eqref{est9} and \eqref{est10} from~\eqref{est4} in the proof 
of Lemma~\ref{l1}, we deduce from this the two estimates
\begin{eqnarray}
\|\partial_s \xi\|_{L^2(I_{T'},H_0)} &\leq& \kappa \bigl( \|\partial_s \beta\|_\infty + c_1' \bigr) , 
                                                                            \label{sim3} \\
\|\xi\|_{L^2(I_{T'},H_1)} &\leq& \kappa \bigl( \|\partial_s \beta\|_\infty + c_1' \bigr) . 
                                                                            \label{sim4}
\end{eqnarray}
Since $\xi=\partial_s w$ we can combine \eqref{sim3} with \eqref{esti} to
\begin{equation}\label{sim5}
\|w\|_{W^{2,2}(I_{T'},H_0)} \,\leq\, \kappa \sqrt{\bigl( \|\partial_s \beta\|_\infty+c_1' \bigr)^2+1}
\end{equation}
and \eqref{sim4} with \eqref{esti} to
\begin{equation}\label{sim6}
\|w\|_{W^{1,2}(I_{T'},H_1)} \,\leq\, \kappa \sqrt{\bigl( \|\partial_s \beta\|_\infty+c_1' \bigr)^2+1}.
\end{equation}
Finally, using $(\cV 3')$, the equation $\pp_s w = \cV(w)$ and also~\eqref{sim6},
we estimate
\begin{eqnarray} \label{sim7}
\|w\|_{L^2(I_{T'},H_2)} &\leq& 
        c_1' \left( \|\partial_s w\|_{L^2(I_{T'},H_1)} + \|w\|_{L^2(I_{T'},H_1)}+1 \right) \\ 
&\leq& c_1' \left( 2 \| w\|_{W^{1,2}(I_{T'},H_1)} +1 \right) \nonumber \\ 
&\leq& c_1' \left( 2 \kappa\sqrt{\bigl(\|\partial_s \beta\|_\infty+c_1' \bigr)^2 +1} +1 \right). \nonumber
\end{eqnarray}
The estimates \eqref{sim5}, \eqref{sim6}, and \eqref{sim7} imply the lemma. 
\proofend

\ni
\textbf{Proof of Proposition~\ref{pr}. } 
By Lemma~\ref{l2} and Lemma~\ref{l3}, for every $0<T'<T$ there
exists a constant $c=c(\kappa,T')$ such that
$$\|w\|_{\bigcap_{k=0}^2 W^{k,2}(I_{T'},H_{2-k})} \leq c.$$
Since by Lemma~\ref{l4} the Hilbert space $\bigcap_{k=0}^2 W^{k,2}(I_{T'},H_{2-k})$ 
compactly and hence continuously embeds into 
the Banach space $\bigcap_{k=0}^1 C^k(I_{T'},H_{1-k})$, Proposition~\ref{pr} follows. 
\proofend

\ni
\textbf{Proof of Theorem~\ref{strong}. }
Theorem~\ref{strong} follows from Theorem~\ref{th} with the help of Proposition~\ref{pr}, 
or alternatively in the same way as Theorem~\ref{th} by combining Lemma~\ref{l2} and
Lemma~\ref{l3} with Lemma~\ref{l4}. 
\proofend

\section{The case of Floer's equation} \label{s:floer}

In this section we explain how classical and less classical Hamiltonian 
and Lagrangian Floer theory fits into the results of this paper
in a special case:
We consider the Floer equation 
on $\R^{2n}$ with the standard symplectic structure $\omega_0$
and a smooth $\omega_0$-compatible almost complex structure~$J$.
This means that for each $p \in \R^{2n}$, $J(p) \colon \R^{2n} \to \R^{2n}$
is a complex structure ($J(p)^2=-\id$) such that 
$\omega_p (\cdot, J(p) \2 \cdot)$ is an inner product on $\R^{2n}$.
The Floer equation then reads
\begin{equation} \label{e:floer}
\pp_s w + J(w) \pp_t w  + \nabla H_t (w) =0
\end{equation}
where $w(s,t)$ is a map from the cylinder $\R \times S^1$ to~$\R^{2n}$, 
and where $H \colon \R^{2n} \times S^1 \to \R$ is a smooth function
and $\nabla$ is the gradient with respect to the Riemannian metric 
$\omega (\cdot, J \cdot)$.

\begin{rem}
{\rm
This situation is enough to obtain compactness for solutions to Floer's equation 
in a general symplectic manifold $(M,\omega)$ for ``short loops'':
If $J_M$ is an $\omega$-compatible almost complex structure and
$\phi \colon (U,\omega_0) \to (M,\omega)$ is a symplectic chart, 
then Floer's equation near a loop $x \subset \phi (U)$ in the chart~$U$ is
the above equation with $J=\phi^*J_M$.
\diam
}
\end{rem}

%

\subsection{Hamiltonian Floer homology on $\R^{2n}$} \label{ss:Ham}

For further use we consider the more general equation
\begin{equation} \label{e:floerX}
\pp_s w + J(w) \pp_t w  + \cX_t (w) =0
\end{equation}
where $\cX \colon  \R^{2n}  \times S^1\to \R^{2n}$ is a smooth vector field.
Choose a smooth map $\Psi \colon \R^{2n} \to \GL (\R^{2n})$
into the space of invertible $2n$-matrices 
such that
\begin{equation} \label{e:Psi}
i \circ \Psi (p) \,=\, \Psi (p) \circ J(p) \quad \forall \, p \in \R^{2n}. 
\end{equation}
We identify $\R^{2n}$ with $\C^n$ and
look at the Sobolev spaces 
$$
W^{k,2}(S^1,\C^n) \,=\, 
\left\{ x(t) = \sum_{j \in \Z} e^{2\pi j t i} \, x_j \;\Big|\,
\sum_{j \in \Z} j^{2k}\,|x_j|^2 < \infty
\right\} , \quad k \in \Z.
$$
Here, $i$ stands for the matrix of the usual complex structure $i \oplus \dots \oplus i$ on $\C^n \cong \R^{2n}$,
and $x_j \in \R^{2n}$. For $k\in \N_0$ we consider inner products on $W^{k,2}$ defined by
$$
\langle x,x \rangle_{W^{k,2}} \,=\, \sum_{j \in \Z} (2\pi j + \tfrac 12)^{2k} \, |x_j|^2 .
$$
The inner product on $W^{0,2}=L^2$ is the usual one. These norms make the operator
$$
\cF (x) \,:=\, -i \2 \pp_t x + \tfrac 12 x
$$
an isometry $W^{k+1,2} \to W^{k,2}$. For convenience, we shall work with these inner products below. The induced norm on $W^{k,2}$ is equivalent to the usual norm given by 
$\|x\|^2_{W^{k,2}} = \sum_{j=0}^k \|\pp_t^j x\|^2_{L^2}$, $k \in \N_0$. In the sequel we shall often use that $W^{k,2}(S^1,\R)$, $k\geq1$, is a Banach algebra, 
in the sense that $\| f g \| \leqslant c \2 \| f \| \|g\|$
for a universal constant~$c$. 

\m \ni 
{\bf The moving frame $\Phi$.}
We first set $H_k = W^{k+1,2}(S^1,\C^n)$. 
For instance, $H_1 = W^{2,2}$ and $H_0 = W^{1,2}$.
For $x \in H_1$ and $v \in H_0$ define the linear map $\Phi (x) \colon H_0 \to H_0$ by
\begin{equation} \label{e:Pmf}
\bigl( \Phi (x) v \bigr)(t) \,=\, \Psi (x(t)) \2 v(t).
\end{equation}
Note that for $x \in H_1$ the coefficients of the $2n$-matrix $\Psi (x)$ are in $H_1 \subset H_0$.
Since $W^{1,2}(S^1,\R)$ is a Banach algebra, 
$\Phi (x)$ thus indeed takes values in~$H_0$.
Similarly, for all $k \ge 1$, \eqref{e:Pmf} defines a map $x \mapsto \Phi (x)$ from~$H_k$ to the linear maps $H_{k-1} \to H_{k-1}$.

\begin{prop} \label{p:mf}
$\Phi$ is a moving frame.
\end{prop}

\proof
$(\Phi 1)$: Fix $x \in H_1$. The map $\Phi (x) \colon H_0 \to H_0$, $v \mapsto \Psi (x(t)) \1 v(t)$
is linear. We shall show that it is bounded. Its inverse $v \mapsto \Psi (x(t))^{-1} v(t)$
is then also bounded (by the bounded inverse theorem, or by the same argument).

For $t \in S^1$ set $A(t) := \Psi (x(t)) \in \GL (\R^{2n})$.
Since $x \in H_1 \subset C^1(S^1,\R^{2n})$ is bounded, there exists a constant~$c$ such that
$\|A(t)\| \leq c$ and $\| \dot A(t)\| \leq c$ for all $t \in S^1$.
Hence
\begin{eqnarray}\label{e:3c}\nonumber
\| \Phi (x) v \|^2_{H_0} &=& \nonumber
\int_{S^1} \|A(t) v(t)\|^2 \, dt + \int_{S^1} \|\dot A(t) v(t) + A(t) \dot v (t) \|^2 \, dt\\\nonumber
&\leq& 3 c^2 \int_{S^1} \left( \|v(t)\|^2 + \|\dot v(t) \|^2 \right) dt \\
&=& 3c^2 \2\|v \|^2_{H_0} .
\end{eqnarray}
A similar computation shows that the map $x \mapsto \Phi (x)$ takes $H_k$ to $\mathcal{L}(H_k,H_k)$
for every $k \geq 0$.
Properties $(\Phi 2)$--$(\Phi 5)$ will readily follow from 

\begin{lemma} \label{le:smooth}
$\Phi \in C^\infty (H_k, \mathcal{L} (H_k,H_k))$ for every $k \geq 0$.
\end{lemma}

\proof
Let $x,\hat x \in H_k$. Since $\Psi \colon \R^{2n} \to \GL(\R^{2n})$ is smooth, for every $t \in S^1$
there exists $\vartheta (t) \in (0,1)$ such that
$$
\Psi \left( x(t) +\hat x (t) \right) \,=\, \Psi (x(t)) + D \Psi \bigl( x(t)+\vartheta (t) \hat x(t) \bigr) \2 \hat x (t)
$$
by the mean value theorem.
With $\hat x_\vartheta (t) = \vartheta(t) \hat x (t)$ we therefore have
\begin{eqnarray} \label{e:Psixhat}
\Phi (x+\hat x) - \Phi (x) - D\Phi (x) (\hat x) &=&
\Psi (x+\hat x) - \Psi (x) - D\Psi (x) (\hat x) \notag \\
&=& \bigl( D\Psi (x+\hat x_\vartheta) - D\Psi (x) \bigr) (\hat x) . 
\end{eqnarray}
Since the inclusion $H_k \subset C^0 := C^0(S^1,\R^{2n})$ is continuous, 
there exists a constant $c$ such that
$\|y\|_{C^0} \leq c \2 \|y\|_{H_k}$ for all $y \in H_k$.
Hence $\|x\|_{C^0} \leq c \2 \|x\|_{H_k}$ and 
$$
\| (x+\hat x_\vartheta)(t) - x(t)\| \,\leq\, \|\hat x_\vartheta \|_{L^\infty}\,\leq\,\|\hat x \|_{L^\infty}= \|\hat x \|_{C^0} 
\,\leq\, c \2 \|\hat x \|_{H_k} \quad \mbox{for all}\, t \in S^1.
$$  
Since $\Psi \colon \R^{2n} \to \GL(\R^{2n})$ is smooth, it follows from \eqref{e:Psixhat} that
$$
\left\| \Phi (x+\hat x) - \Phi (x) - D\Phi (x) (\hat x) \right\|_{\mathcal{L}(H_k,H_k)}\,=\, 
o(\|\hat x\|_{H_k}) .
$$
We have shown that $\Phi \colon H_k \to \mathcal{L}(H_k,H_k)$ is differentiable for every $k\in \N$. 
Since $\Psi$ is smooth, we can iterate this argument and find that
$\Phi \in C^\infty (H_k, \mathcal{L}(H_k,H_k))$ for every $k\in \N$.
\proofend

The map in $(\Phi 2)$ is the composition 
$$
H_1 \times H_0 \,\subset\, H_0 \times H_0 \,\to\, \mathcal{L}(H_0,H_0) \times H_0 \,\to\, H_0 
$$
where the second map is $(x,v) \mapsto (\Phi (x),v)$ is well-defined by Lemma~\eqref{le:smooth} and the third map is $(A,v) \mapsto Av$.
The inclusion and the third map are continuous, and the second map is continuous by Lemma~\ref{le:smooth}.
In fact, by this lemma, $\Phi \in C^1(H_0, \mathcal{L}(H_0,H_0))$,
which in particular implies
properties $(\Phi 3)$ and~$(\Phi 4)$.
Property~$(\Phi 5)$ also follows from Lemma~\ref{le:smooth}.
Finally,    
$(\Phi 6)$ holds because for $x \in H_1$ with $\|x\|_1 \leq \kappa$ we have $\|x\|_{C^0} \leq c \kappa$
and since the coefficients of the matrices $\Psi (p)$, $\Psi(p)^{-1}$ together with their derivatives up to order two 
are uniformly bounded on the ball $\{ p \in \C^n : \|p\| \leq c \kappa\}$.
By means of example we spell this out this for $\| \Phi(x)v\|_{H_1} \leq c_0\2 \|v\|_{H_1}$ 
for all $x \in H_1$ with $\|x\|_1 \leq \kappa$ and $v \in H_1$:
Using $\Phi(x)v = \Psi (x) v$ we have
\begin{eqnarray*}
\|\Phi (x) v\|_{H_1}^2 &=&  \|\Psi (x) v\|_{L^2}^2 + \|D \Psi (x) (\dot x, v) + \Psi(x)\dot v \|_{L^2}^2 + \\
&& \|D^2 \Psi (x) (\dot x, \dot x, v) + D \Psi(x) (\ddot x, \dot v) + 2 D\Psi(x)(\dot x, \dot v) + \Psi (x) \ddot v \|_{L^2}^2 \\
&\leq& 
\|\Psi (x) v\|_{L^2}^2 + 2 \left( \|D \Psi (x) (\dot x, v) \|_{L^2}^2 + \| \Psi(x) \dot v \|_{L^2}^2 \right) + \\
&& 8 \left( \|D^2 \Psi (x) (\dot x, \dot x, v) \|_{L^2}^2 + \| D \Psi(x) (\ddot x, \dot v) \|_{L^2}^2 + \| D\Psi(x)(\dot x, \dot v) \|_{L^2}^2 + \| \Psi (x) \ddot v \|_{L^2}^2 \right) 
\end{eqnarray*}
Integrating by parts we have 
$\int_{S^1} \langle \ddot x, \dot v \rangle \2 dt = -\int_{S^1} \langle \dot x, \ddot v\rangle \2 dt$.
Since also $\|x\|_{C^1} \leq c \2 \|x\|_1 \leq c \kappa$ for a universal constant~$c$,
we find a constant $c_0$ depending only on $\kappa$ and on~$\Psi$ and its derivatives up to order two 
such that the above some is bounded by
$c_0^2 \left( \|v\|_{L^2}^2 + \| \dot v\|_{L^2}^2 + \| \ddot v\|_{L^2}^2 \right) = c_0^2 \2 \|v\|_{H_1}^2
$. 
\proofend

\m \ni 
{\bf The unregularized vector field $\cV$.}
In view of~\eqref{e:floerX} we define the vector field $\cV \colon H_1 \to H_0$ by 
\begin{equation} \label{e:uv}
\cV (x) \,=\, -J(x) \2 \pp_t x - \cX_t (x) .
\end{equation}
Since $J \colon \R^{2n} \to \End (\R^{2n})$ and $\cX \colon S^1 \times \R^{2n} \to \R^{2n}$
are smooth and since $W^{1,2}(S^1,\R^{2n})$ is a Banach algebra, $\cV$ indeed takes values in~$H_0$,
and it readily follows that $\cV \in C^0(H_1,H_0) \cap C^0(H_2,H_1)$.

The differential $D \cV (x) \colon H_1 \to H_0$ at $x \in H_1$ is given by
\begin{equation} \label{e:DcV}
D \cV (x) (\hat x) \,=\, -DJ (x) \left( \hat x, \pp_t x \right) - J(x) \pp_t \hat x
- D \cX_t (x) (\hat x) .
\end{equation}
From this, the continuity of the map $H_1 \times H_1 \to H_0$, $(x,\hat x) \mapsto D\cV(x)\hat x$
asked in~$(\cV 1)$ is immediate.
To verify $(\cV 2)$ we compute, using~\eqref{e:DcV},
\begin{eqnarray*}
\bigl( \Phi (x) D\cV (x) \Phi (x)^{-1} (\hat x) \bigr) (t) &=&
- \Psi (x(t)) \2  DJ (x(t)) \bigl( \Psi (x(t))^{-1} \hat x (t), \pp_t x(t) \bigr) \\
&&
- \Psi (x(t)) \2 J (x(t)) \, \pp_t \bigl( \Psi (x(t))^{-1} \hat x (t) \bigr) \\
&&
- \Psi (x(t)) D \cX_t (x(t)) \, \Psi (x(t))^{-1} \hat x (t) .
\end{eqnarray*}
Since 
\begin{eqnarray*}
\pp_t \bigl( \Psi (x(t))^{-1} \hat x (t) \bigr) &=& 
-\Psi (x(t))^{-1} D \Psi (x(t)) \bigl( \pp_t x(t), \Psi (x(t))^{-1} \hat x (t) \bigr)
\\
&& + \Psi (x(t))^{-1} \pp_t \hat x (t)
\end{eqnarray*}
and by the defining property~\eqref{e:Psi} of $\Psi$,
\begin{eqnarray*}
\Psi (x(t)) \2 J(x(t)) \, \pp_t \bigl( \Psi (x(t))^{-1} \hat x (t) \bigr)
&=&
- i \2 D \Psi (x(t)) \bigl( \pp_t x(t), \Psi (x(t))^{-1} \hat x (t) \bigr) \\
&&
+ i \2 \pp_t \hat x (t) .
\end{eqnarray*}
Hence
\begin{eqnarray*}
\bigl( \Phi (x) D\cV (x) \Phi (x)^{-1} (\hat x) \bigr) (t) &=&
- \Psi (x(t)) \2  DJ (x(t)) \bigl( \Psi (x(t))^{-1} \hat x (t), \pp_t x(t) \bigr) \\
&&
+ i \2 D \Psi (x(t)) \bigl( \pp_t x(t), \Psi (x(t))^{-1} \hat x (t) \bigr) \\
&&
- i \2 \pp_t \hat x (t)  \\
&&
- \Psi (x(t)) D \cX_t (x(t)) \, \Psi (x(t))^{-1} \hat x (t) .
\end{eqnarray*}
It is tempting to take $\cF (\hat x) = -i \2 \pp_t \hat x$. 
However, this Fredholm operator $H_1 \to H_0$ of index zero has a 1-dimensional kernel. 
Since we need $\cF$ to be invertible, we define
$$
\cF (\hat x) \,:=\, -i \2 \pp_t \hat x + \tfrac 12 \hat x .
$$
Then 
\begin{eqnarray} \label{e:Pdef}
\bigl( \Phi (x) D\cV (x) \Phi (x)^{-1} - \cF \bigr)  (\hat x) (t) 
&=&
-
\Psi (x(t)) \2  DJ (x(t)) \bigl( \Psi (x(t))^{-1} \hat x (t), \pp_t x(t) \bigr) \notag \\
&&
+ i \2 D \Psi (x(t)) \bigl( \pp_t x(t), \Psi (x(t))^{-1} \hat x (t) \bigr) \notag \\
&&
- \tfrac 12 \2 \hat x (t)  \notag \\
&&
- \Psi (x(t)) D \cX_t (x(t)) \, \Psi (x(t))^{-1} \hat x (t) .
\end{eqnarray}
We can thus define $\cP (x) \colon H_0 \to H_0$ for $x \in H_1$ by \eqref{e:Pdef},
and then $(\cV 2)$ holds true.
(Use again the Banach algebra structure of $W^{1,2}(S^1,\R)$ to see that $\cP (x)$ takes values in~$H_0$
and that $H_1 \times H_0 \to H_0$, $(x,\hat x) \mapsto \cP(x)\hat x$ is continuous.) 

To verify $(\cV 3)$ fix $\kappa >0$.
The existence of a constant $c_1(\kappa)$ such that 
$$
\|\cP (x) \hat x\|_{H_0} \,\leq\, c_1(\kappa)
$$
for all $x \in H_1$ with $\|x\|_{H_1} \leq \kappa$ and all $\hat x \in H_0$ with $\|\hat x\|_{H_0} \leq 1$
follows from the smoothness of $\Psi$ and~$J$ and from the Banach algebra structure of $W^{1,2}(S^1,\R)$.

Further, since $\|x\|_{C^0} \le c\2 \|x\|_{H_1} \leq c \2 \kappa$ and
since $J^{-1}=-J$ and $\cX_t$ are smooth, there are constants $d_1=d_1(\kappa)$ and $d_2 = d_2(\kappa)$
such that for all $x \in H_2$ with $\|x\|_{H_1} \leq \kappa$:
\begin{eqnarray*}
\|\cV(x)\|_{H_1} &=& \| J(x) \partial_t x + \cX_t(x) \|_{H_1} \\
&\geq& \| J(x) \partial_t x \|_{H_1} - \| \cX_t(x) \|_{H_1} \\
&\geq& d_1 \|x\|_{H_2} - d_2,
\end{eqnarray*}
and so $\|x\|_{H_2} \leq \frac{1}{d_1} \|\cV(x) \|_{H_1} + \frac{d_2}{d_1}$.

\b \ni 
{\bf The case of elementary unregularized vector fields $\cV$.}
We next look at the case of elementary vector fields. 
That the moving frame defined by~\eqref{e:Pmf} is the identity means that $\Psi \equiv \id$,
that is, $J \equiv i$. Hence
$$
\cV (x) \,=\, -i \pp_t x - \cX_t(x) 
$$
with $D \cV (x)\hat x = -i \pp_t \hat x - D\cX_t (x)(\hat x)$. Then 
$$
\cP (x)\hat x \,=\, D\cV (x) \hat x - \cF \hat x \,=\, - \tfrac 12 \hat x - D\cX_t (x)(\hat x).
$$
Since all non-linear terms of $\cV$, $D\cV$ and~$\cP$ 
are now in the summand involving the vector field~$\cX$,
that contains no derivatives of~$x$ or~$\hat x$,
we can now take $H_k = W^{k,2}(S^1,\C^n)$,
so that $H_0 = L^2$, $H_1 = W^{1,2}$, $H_2 = W^{2,2}$.
Properties~$(\cV 1)$ and~$(\cV 2)$ are readily checked.
(Note that the multiplication of an $L^2$ function with a $W^{1,2}$ function lies in~$L^2$.)

To verify $(\cV 3')$, we need the following assumption on~$\cX$:
There exists $\gamma \in \R \setminus \{0\}$ and a constant $c=c(\cX)>0$ such that 
for all $t \in S^1$ and all $p \in \C^n$,
\begin{equation} \tag{$\cX$}
\| \pp_t \cX_t (p)\| \leq c,
\quad
\| \cX_t(p) - \gamma p \| \leq c,
\quad 
\| D \cX_t (p) \| \leq c .
\end{equation}

By the third point, $\|\cP(x)\|_{\mathcal{L}(H_0,H_0)} \leq \frac 12 + c$.
By the second point in~$(\cX)$ we can estimate 
\begin{eqnarray*}
\|\cV (x)\|_{L^2} &=& \| i\2 \pp_t x + \cX_t (x) \|_{L^2} \\
&\geq& \| i\2 \pp_t x + \gamma x \|_{L^2} - \| \cX_t(x)-\gamma x \|_{L^2} \\
&\geq&
\min \{1,|\gamma|\} \, \|x\|_{W^{1,2}} -c .
\end{eqnarray*}
Finally, using all three assumptions in $(\cX)$,
\begin{eqnarray*}
\|\cX_t(x)\|_{W^{1,2}} &\leq& \|\cX_t(x)\|_{L^2} + \| \pp_t \bigl( \cX_t(x) \bigr) \|_{L^2} \\
&\leq& \|\cX_t(x)\|_{L^2} + \| (\pp_t \cX) (x) \|_{L^2} + \| D \cX_t (x) (\pp_t x) \|_{L^2} \\
&\leq& c+\gamma \|x\|_{L^2} + c + c \2 \|\pp_tx\|_{L^2} \\
&\leq& 2c+\gamma \|x\|_{L^2} + c \2 \|x\|_{W^{1,2}}.
\end{eqnarray*}
Therefore, 
\begin{eqnarray*}
\|\cV (x)\|_{W^{1,2}} &\geq& \| i\2 \pp_t x + \gamma x \|_{W^{1,2}} - \| \cX_t (x) - \gamma x \|_{W^{1,2}} \\
&\geq&
\min \{1,|\gamma|\} \, \|x\|_{W^{2,2}} - 2c - (c+2\gamma) \|x\|_{W^{1,2}}.
\end{eqnarray*}
Property $(\cV 3')$ therefore holds with 
$$
c_1' \,:=\ \frac{\max\{ 1,2c,c+2|\gamma| \}}{\min \{ 1,|\gamma| \}} .
$$

We have verified that Theorems~\ref{th} and~\ref{strong}
hold true for vector fields~$\cV$ of the form~\eqref{e:floerX},
under the assumption on the bounds~\eqref{e:ass1} and~\eqref{e:ass2}.
In the Hamiltonian case $\cX_t = \nabla H_t$, assumption~$(\cX)$
becomes 
\begin{equation} \tag{H}
\|\nabla \pp_t H_t (p)\| \leq c,
\quad
\| \nabla H_t(p) - \gamma p \| \leq c,
\quad 
\| \Hess H_t (p) \| \leq c .
\end{equation}
For instance, 
functions of the form $H_t(p) = \frac{\gamma}{2} \2 \|p\|^2 + ap + c + f(t,p)$ where $f$ has compact support satisfy~(H).
Recall from Remark~\ref{re:bound} that in the Hamiltonian case the bounds~\eqref{e:ass1} and \eqref{e:ass2} are known.

\begin{rem} \label{re:t}
{\rm
In Floer homology one often makes the compatible almost complex structure~$J$
depend on~$t$, in order to achieve transversality. 
The previous arguments go through for such families~$J_t$ 
(use families $\Psi_t(p)$ conjugating $J_t(p)$ to~$i$). 
However, to achieve transversality in our setting 
it will be more convenient and more natural to perturb the whole vector field~$\cV$. 
}
\end{rem}

\subsection{Lagrangian boundary conditions} \label{ss:Lag}
Our framework is general enough to embrace many classical Floer theories. 
We discuss one more example.
The local model for Lagrangian Floer homology is the pair $(\R^{2n},\R^n)$
where $\R^n$ is the real part of $\C^n = \R^n \oplus i \R^n$.
The Floer equation is again~\eqref{e:floer}, where now $w(s,t)$ is a map from the strip 
$\R \times [0,1]$ to~$\R^{2n}$ mapping the boundary lines to~$\R^n$.
We thus take for $k \in \Z$ the Hilbert spaces
$$
W^{k,2}([0,1],\R^{2n},\R^n) \,:=\, 
\left\{ x(t) = \sum_{\ell \in \Z} e^{\pi j t i} x_j \,\Big|\,
t \in [0,1],\; x_j \in \R^n,\; \sum_{j \in \Z} |j|^{2k} \, \|x_j\|^2 < \infty \right\} 
$$
of paths with endpoints on~$\R^n$. 
Define~$\Phi$ and~$\cV$ again by \eqref{e:Pmf} and \eqref{e:uv}, with $\cX_t = \nabla H_t$,
and take again $H_k = W^{k+1,2}$ in the general case and $H_k = W^{k,2}$ in the case that $\Phi$ is trivial.
Arguing literally as before we then see that $\Phi$ is a moving frame and that $\cV$ is an unregularized vector field, that in the case where $\Phi$ is trivial is elementary
if $H$ satisfies the growth condition~(H).
The compactness theorems~\ref{th} and~\ref{strong} thus also hold true for the Lagrangian Floer equation.

The following lemma describes the elements of $W^{k,2}([0,1],\R^{2n},\R^n)$ for $k \geq 1$
in a more geometric way. 
For $k \geq 1$
define $W^{k,2}_{\mathrm{bc}}([0,1],\R^{2n})$ as the space of paths 
$x \in W^{k,2}([0,1],\R^{2n})$ for which
\begin{eqnarray} 
\pp_t^\ell x(0),\, \pp_t^\ell x(1) \in \phantom{i} \R^n & 
                                \mbox{ if $\ell$ is even,} & 0 \leq \ell \leq k-1,  \label{e:even}\\
\pp_t^\ell x(0),\, \pp_t^\ell x(1) \in i\R^n            & 
                                \mbox{ if $\ell$ is odd, \,}& 0 \leq \ell \leq k-1. \label{e:odd}
\end{eqnarray}
These spaces where introduced 
by Tatjana Simcevic~\cite{Sim14} in her Hardy space approach to Lagrangian Floer gluing.

\begin{lemma}
For $k \geq 1$,
$W^{k,2}([0,1],\R^{2n},\R^n) = W^{k,2}_{\mathrm{bc}}([0,1],\R^{2n})$. 
\end{lemma}

\proof
For $x = \sum_{\ell \in \Z} e^{\pi j t i} x_j \in W^{k,2}([0,1],\R^{2n},\R^n)$ and $\ell \in \{0,\dots,k-1\}$
we have
$$
\pp_t^\ell x (t) \,=\, i^\ell \sum_{j \in \Z} (\pi j)^\ell e^{\pi jti} x_j .
$$
Since $e^{\pi j i} = (-1)^j \id$, \eqref{e:even} and \eqref{e:odd} hold true.

Now assume that $x \in W^{k,2}([0,1],\R^{2n})$ satisfies \eqref{e:even} and \eqref{e:odd}.
Let $S^1(2) = \R / 2 \Z$ be the circle of length~$2$, and let $\gamma_x \colon S^1(2) \to \R^{2n}$ 
be the loop obtained by reflecting $x$ at~$\R^n$:
\begin{eqnarray} 
\gamma_x(t) &=& x(t)                \qquad \;\;\, \mbox{ if }\, t \in [0,1],  \label{e:ga0}\\
\gamma_x(t) &=&  \overline{x(2-t)}  \quad \mbox{ if }\, t \in [1,2] . \label{e:ga1}
\end{eqnarray}
Then 
\begin{equation} \label{e:gabar}
\gamma_x (2-t) \,=\, \overline{\gamma_x(t)} \quad \mbox{ for all }\, t \in [0,2] .
\end{equation}
We claim that $\gamma_x \in W^{k,2}(S^1(2), \R^{2n})$.
Indeed, $\gamma_x$ has $k$ weak derivatives in $L^2(S^1(2),\R^{2n})$
since $x$ has $k$ weak derivatives in $L^2([0,1],\R^{2n})$.
Further, $\gamma_x$ has $k-1$ continuous derivatives on $S^1(2) \setminus \{0,1\}$
since $x$ has $k-1$ continuous derivatives on $(0,1)$.
Geometrically, it is clear from~\eqref{e:gabar} and \eqref{e:even}, \eqref{e:odd} that
$\gamma_x$ also has $k-1$ continuous derivatives at $0$ and~$1$.
To see this formally, we use that $x$ has $k-1$ continuous derivatives at $0$ and~$1$
and for $\ell \leq k-1$ compute that at~$t=1$,
%
$$
\lim_{t \nearrow 1} \pp_t^\ell \, \gamma_x (t) 
\stackrel{\eqref{e:ga0}}{=}
\lim_{t \nearrow 1} \pp_t^\ell \, x (t) 
\stackrel{\eqref{e:even},\eqref{e:odd}}{=}  
(-1)^\ell \lim_{t \nearrow 1} \pp_t^\ell \, \overline{x (t)} 
\,=\, \lim_{t \searrow 1} \pp_t^\ell \, \overline{x (2-t)} 
\stackrel{\eqref{e:ga1}}{=}
\lim_{t \searrow 1} \pp_t^\ell \, \gamma_x{t}
$$
and similarly $\pp_t^\ell \, \gamma_x (t)$ is continuous at $t=0$.

Since $\gamma_x \in W^{k,2}(S^1(2),\R^{2n})$ we can write
$$
\gamma_x(t) \,=\, \sum_{j \in \Z} e^{\pi jti} \,x_j, \quad t \in S^1(2),\; x_j \in \C^n.
$$
With $\sum_{j \in \Z} |j|^{2k} \|x_j\|^2 < \infty$.
Property \eqref{e:gabar} then becomes
$\sum_{j \in \Z} e^{\pi j(2-t)i} \,x_j = \sum_{j \in \Z} e^{-\pi jti} \,\overline{x_j}$.
Since $e^{\pi j 2 i} = \id$ we find that $x_j \in \R^n$ for all~$j$.
\proofend

\subsection{Delay equations}  \label{ss:delay}
A delay equation is a differential equation in which the velocity does not only depend on the
present state but also on states in the past. 
The simplest case on~$\R^{2n}$ is the differential equation   
$$
\dot x (t) \,=\ \sum_{j=1}^m \cX^j_t (x(t-\tau_j)) 
$$
where $\cX^j \colon S^1 \times \R^{2n} \to \R^{2n}$ are vector fields and 
$0 \leq \tau_1 < \tau_2 < \dots < \tau_m$.
Since the curves $x(t-\tau_j)$ have the same norms in $W^{k,2}(S^1,\R^{2n})$ and $C^k(S^1,\R^{2n})$
as~$x(t)$, the arguments in~\S \ref{ss:Ham} show that
$$
\cV (x)(t) \,:=\, - J(x(t)) \pp_t x(t) - \sum_{j=1}^m \cX^j_t (x(t-\tau_j)) 
$$
is an (elementary) unregularized vector field.
For Hamiltonian delay vector fields on~$\R^{2n}$, for which also assumptions~\eqref{e:ass1}
and~\eqref{e:ass2} are verified, we refer to~\cite{albers-frauenfelder-schlenk-delay}.

\end{document}